\pdfoutput=1

\documentclass[12pt,a4paper]{article}

\usepackage{amsfonts,amssymb,amsthm,amsmath,latexsym}
\usepackage{eqsection,indent}
\usepackage{subeqnarray, eepicemu, url,cite,bm}
\usepackage{multirow}  % See http://en.wikibooks.org/wiki/LaTeX/Tables
\usepackage{tensor}  % To help do 2F1, 2\phi1, etc
\usepackage{graphicx,graphbox}
\usepackage{float}  % http://tex.stackexchange.com/questions/2275/keeping-tables-figures-close-to-where-they-are-mentioned

\date{March 2, 2021 \\ revised July 27, 2021 \\[1mm]
   To appear in the {\em Proceedings of the AMS}}   %% REMEMBER TO PUT FINAL DATE HERE!!!!!
\begin{document}

\title{Multiple Laguerre polynomials: \\[1mm]
\hspace*{-0.5cm}\hbox{Combinatorial model and Stieltjes moment representation}}

\author{
     {\small Alan D.~Sokal}                  \\[2mm]
     {\small\it Department of Mathematics}   \\[-2mm]
     {\small\it University College London}   \\[-2mm]
     {\small\it Gower Street}                \\[-2mm]
     {\small\it London WC1E 6BT}             \\[-2mm]
     {\small\it UNITED KINGDOM}              \\[-2mm]
     {\small\tt sokal@math.ucl.ac.uk}        \\[-2mm]
     {\protect\makebox[5in]{\quad}}  % To force authors' names to be written
                                     %   vertically, one above another.
                                     % (\author seems to put them side-by-side
                                     %   if there is room.)
     \\[-2mm]
     {\small\it Department of Physics}       \\[-2mm]
     {\small\it New York University}         \\[-2mm]
     {\small\it 726 Broadway}                \\[-2mm]
     {\small\it New York, NY 10003}          \\[-2mm]
     {\small\it USA}                         \\[-2mm]
     {\small\tt sokal@nyu.edu}               \\[3mm]
}

\maketitle
\thispagestyle{empty}   % Suppress page number on front page.

\begin{abstract}
I give a combinatorial interpretation of the multiple Laguerre polynomials
of the first kind of type~II, generalizing the digraph model found by
Foata and Strehl for the ordinary Laguerre polynomials.
I also give an explicit integral representation for these polynomials,
which shows that they form a multidimensional Stieltjes moment sequence
whenever $x \le 0$.
\end{abstract}


\bigskip
\bigskip
\noindent
{\bf Mathematics Subject Classification (MSC 2010) codes:}
33C45 (Primary); 05A15, 05A19, 30E05, 42C05 (Secondary).

%%\clearpage
\vspace*{1cm}

\newtheorem{theorem}{Theorem}[section]
\newtheorem{proposition}[theorem]{Proposition}
\newtheorem{lemma}[theorem]{Lemma}
\newtheorem{corollary}[theorem]{Corollary}
\newtheorem{definition}[theorem]{Definition}
\newtheorem{conjecture}[theorem]{Conjecture}
\newtheorem{question}[theorem]{Question}
\newtheorem{problem}[theorem]{Problem}
\newtheorem{openproblem}[theorem]{Open Problem}
\newtheorem{example}[theorem]{Example}
\newtheorem{remark}[theorem]{Remark}

\renewcommand{\theenumi}{\alph{enumi}}
\renewcommand{\labelenumi}{(\theenumi)}
\def\eop{\hbox{\kern1pt\vrule height6pt width4pt
depth1pt\kern1pt}\medskip}
\def\prf{\par\noindent{\bf Proof.\enspace}\rm}
\def\rmk{\par\medskip\noindent{\bf Remark\enspace}\rm}

\newcommand{\textbfit}[1]{\textbf{\textit{#1}}}

\newcommand{\bigdash}{%
\smallskip\begin{center} \rule{5cm}{0.1mm} \end{center}\smallskip}

% DO NOT USE \par WITHIN A \footnote;  USE \safepar INSTEAD.
\newcommand{\safepar}{ {\protect\hfill\protect\break\hspace*{5mm}} }

\newcommand{\be}{\begin{equation}}
\newcommand{\ee}{\end{equation}}
\newcommand{\<}{\langle}
\renewcommand{\>}{\rangle}
\newcommand{\widebar}{\overline}
\def\reff#1{(\protect\ref{#1})}
\def\spose#1{\hbox to 0pt{#1\hss}}
\def\ltapprox{\mathrel{\spose{\lower 3pt\hbox{$\mathchar"218$}}
    \raise 2.0pt\hbox{$\mathchar"13C$}}}
\def\gtapprox{\mathrel{\spose{\lower 3pt\hbox{$\mathchar"218$}}
    \raise 2.0pt\hbox{$\mathchar"13E$}}}
\def\textprime{${}^\prime$}
\def\proof{\par\medskip\noindent{\sc Proof.\ }}
\def\firstproof{\par\medskip\noindent{\sc First Proof.\ }}
\def\secondproof{\par\medskip\noindent{\sc Second Proof.\ }}
\def\alternateproof{\par\medskip\noindent{\sc Alternate Proof.\ }}
\def\algebraicproof{\par\medskip\noindent{\sc Algebraic Proof.\ }}
\def\graphicalproof{\par\medskip\noindent{\sc Graphical Proof.\ }}
\def\combinatorialproof{\par\medskip\noindent{\sc Combinatorial Proof.\ }}
\def\proofof#1{\bigskip\noindent{\sc Proof of #1.\ }}
\def\firstproofof#1{\bigskip\noindent{\sc First Proof of #1.\ }}
\def\secondproofof#1{\bigskip\noindent{\sc Second Proof of #1.\ }}
\def\thirdproofof#1{\bigskip\noindent{\sc Third Proof of #1.\ }}
\def\algebraicproofof#1{\bigskip\noindent{\sc Algebraic Proof of #1.\ }}
\def\combinatorialproofof#1{\bigskip\noindent{\sc Combinatorial Proof of #1.\ }}
\def\sketchofproof{\par\medskip\noindent{\sc Sketch of proof.\ }}
\renewcommand{\qed}{ $\square$ \bigskip}
\newcommand{\myendremark}{ $\blacksquare$ \bigskip}
\def\half{ {1 \over 2} }
\def\third{ {1 \over 3} }
\def\twothird{ {2 \over 3} }
\def\smfrac#1#2{{\textstyle{#1\over #2}}}
\def\smhalf{ {\smfrac{1}{2}} }
\newcommand{\real}{\mathop{\rm Re}\nolimits}
\renewcommand{\Re}{\mathop{\rm Re}\nolimits}
\newcommand{\imag}{\mathop{\rm Im}\nolimits}
\renewcommand{\Im}{\mathop{\rm Im}\nolimits}
\newcommand{\sgn}{\mathop{\rm sgn}\nolimits}
\newcommand{\tr}{\mathop{\rm tr}\nolimits}
\newcommand{\supp}{\mathop{\rm supp}\nolimits}
\newcommand{\disc}{\mathop{\rm disc}\nolimits}
\newcommand{\diag}{\mathop{\rm diag}\nolimits}
\newcommand{\tridiag}{\mathop{\rm tridiag}\nolimits}
\newcommand{\AZ}{\mathop{\rm AZ}\nolimits}
\newcommand{\NC}{\mathop{\rm NC}\nolimits}
\newcommand{\PF}{{\rm PF}}
\newcommand{\rk}{\mathop{\rm rk}\nolimits}
\newcommand{\perm}{\mathop{\rm perm}\nolimits}
\def\hboxscript#1{ {\hbox{\scriptsize\em #1}} }
\renewcommand{\emptyset}{\varnothing}
\newcommand{\eqdef}{\stackrel{\rm def}{=}}

\newcommand{\restrict}{\upharpoonright}

\newcommand{\compinv}{{\langle -1 \rangle}}   % Compositional inverse

\newcommand{\scra}{{\mathcal{A}}}
\newcommand{\scrb}{{\mathcal{B}}}
\newcommand{\scrc}{{\mathcal{C}}}
\newcommand{\scrd}{{\mathcal{D}}}
\newcommand{\scrdtilde}{{\widetilde{\mathcal{D}}}}
\newcommand{\scre}{{\mathcal{E}}}
\newcommand{\scrf}{{\mathcal{F}}}
\newcommand{\scrg}{{\mathcal{G}}}
\newcommand{\scrh}{{\mathcal{H}}}
\newcommand{\scri}{{\mathcal{I}}}
\newcommand{\scrj}{{\mathcal{J}}}
\newcommand{\scrk}{{\mathcal{K}}}
\newcommand{\scrl}{{\mathcal{L}}}
\newcommand{\scrlbar}{{\overline{\mathcal{L}}}}
\newcommand{\scrm}{{\mathcal{M}}}
\newcommand{\scrn}{{\mathcal{N}}}
\newcommand{\scro}{{\mathcal{O}}}
\newcommand\scroo{
  \mathchoice
    {{\scriptstyle\mathcal{O}}}% \displaystyle
    {{\scriptstyle\mathcal{O}}}% \textstyle
    {{\scriptscriptstyle\mathcal{O}}}% \scriptstyle
    {\scalebox{0.6}{$\scriptscriptstyle\mathcal{O}$}}%\scriptscriptstyle
  }
%% Taken from https://tex.stackexchange.com/questions/191059/how-to-get-a-small-letter-version-of-mathcalo
\newcommand{\scrp}{{\mathcal{P}}}
\newcommand{\scrq}{{\mathcal{Q}}}
\newcommand{\scrr}{{\mathcal{R}}}
\newcommand{\scrs}{{\mathcal{S}}}
\newcommand{\scrt}{{\mathcal{T}}}
\newcommand{\scrv}{{\mathcal{V}}}
\newcommand{\scrw}{{\mathcal{W}}}
\newcommand{\scrz}{{\mathcal{Z}}}
\newcommand{\SP}{{\mathcal{SP}}}
\newcommand{\ST}{{\mathcal{ST}}}

\newcommand{\bfa}{{\mathbf{a}}}
\newcommand{\bfb}{{\mathbf{b}}}
\newcommand{\bfc}{{\mathbf{c}}}
\newcommand{\bfd}{{\mathbf{d}}}
\newcommand{\bfe}{{\mathbf{e}}}
\newcommand{\bfh}{{\mathbf{h}}}
\newcommand{\bfj}{{\mathbf{j}}}
\newcommand{\bfi}{{\mathbf{i}}}
\newcommand{\bfk}{{\mathbf{k}}}
\newcommand{\bfl}{{\mathbf{l}}}
\newcommand{\bfL}{{\mathbf{L}}}
\newcommand{\bfm}{{\mathbf{m}}}
\newcommand{\bfn}{{\mathbf{n}}}
\newcommand{\bfp}{{\mathbf{p}}}
\newcommand{\bfr}{{\mathbf{r}}}
\newcommand{\bft}{{\mathbf{t}}}
\newcommand{\bfu}{{\mathbf{u}}}
\newcommand{\bfv}{{\mathbf{v}}}
\newcommand{\bfw}{{\mathbf{w}}}
\newcommand{\bfx}{{\mathbf{x}}}
\newcommand{\bfX}{{\mathbf{X}}}
\newcommand{\bfy}{{\mathbf{y}}}
\newcommand{\bfz}{{\mathbf{z}}}
\renewcommand{\k}{{\mathbf{k}}}
\newcommand{\n}{{\mathbf{n}}}
\newcommand{\vv}{{\mathbf{v}}}
\newcommand{\bv}{{\mathbf{v}}}
\newcommand{\w}{{\mathbf{w}}}
\newcommand{\x}{{\mathbf{x}}}
\newcommand{\y}{{\mathbf{y}}}
\newcommand{\cc}{{\mathbf{c}}}
\newcommand{\zero}{{\mathbf{0}}}
\newcommand{\one}{{\mathbf{1}}}
\newcommand{\bmm}{{\mathbf{m}}}

\newcommand{\ahat}{{\widehat{a}}}
\newcommand{\Zhat}{{\widehat{Z}}}

\newcommand{\C}{{\mathbb C}}
\newcommand{\D}{{\mathbb D}}
\newcommand{\Z}{{\mathbb Z}}
\newcommand{\N}{{\mathbb N}}
\newcommand{\Q}{{\mathbb Q}}
\newcommand{\PP}{{\mathbb P}}
\newcommand{\R}{{\mathbb R}}
\newcommand{\RR}{{\mathbb R}}
\newcommand{\E}{{\mathbb E}}

\newcommand{\Sym}{{\mathfrak{S}}}
\newcommand{\SymB}{{\mathfrak{B}}}
\newcommand{\Alt}{{\mathrm{Alt}}}

\newcommand{\germanA}{{\mathfrak{A}}}
\newcommand{\germanB}{{\mathfrak{B}}}
\newcommand{\germanQ}{{\mathfrak{Q}}}
\newcommand{\germanh}{{\mathfrak{h}}}

\newcommand{\myle}{\preceq}
\newcommand{\myge}{\succeq}
\newcommand{\mygt}{\succ}

\newcommand{\B}{{\sf B}}
\newcommand{\OB}{B^{\rm ord}}
\newcommand{\OS}{{\sf OS}}
\newcommand{\OO}{{\sf O}}
\newcommand{\OSP}{{\sf OSP}}
\newcommand{\Eu}{{\sf Eu}}
\newcommand{\ERR}{{\sf ERR}}
\newcommand{\sfB}{{\sf B}}
\newcommand{\sfD}{{\sf D}}
\newcommand{\sfE}{{\sf E}}
\newcommand{\sfG}{{\sf G}}
\newcommand{\sfJ}{{\sf J}}
\newcommand{\sfL}{{\sf L}}
\newcommand{\sfLhat}{{\widehat{{\sf L}}}}
\newcommand{\sfLtilde}{{\widetilde{{\sf L}}}}
\newcommand{\sfP}{{\sf P}}
\newcommand{\sfQ}{{\sf Q}}
\newcommand{\sfS}{{\sf S}}
\newcommand{\sfT}{{\sf T}}
\newcommand{\sfW}{{\sf W}}
\newcommand{\sfMV}{{\sf MV}}
\newcommand{\AMV}{{\sf AMV}}
\newcommand{\BM}{{\sf BM}}
\newcommand{\emIB}{B^{\rm irr}}
\newcommand{\emIP}{P^{\rm irr}}
\newcommand{\emOB}{B^{\rm ord}}
\newcommand{\emCB}{B^{\rm cyc}}
\newcommand{\emSC}{P^{\rm cyc}}

\newcommand{\lev}{{\rm lev}}
\newcommand{\stat}{{\rm stat}}
\newcommand{\cyc}{{\rm cyc}}
\newcommand{\e}{{\rm e}}
\newcommand{\ecyc}{{\rm ecyc}}
\newcommand{\epa}{{\rm epa}}
\newcommand{\iv}{{\rm iv}}
\newcommand{\pa}{{\rm pa}}
\newcommand{\mysteryone}{{\rm mys1}}
\newcommand{\mysterytwo}{{\rm mys2}}
\newcommand{\Asc}{{\rm Asc}}
\newcommand{\asc}{{\rm asc}}
\newcommand{\Des}{{\rm Des}}
\newcommand{\des}{{\rm des}}
\newcommand{\Exc}{{\rm Exc}}
\newcommand{\exc}{{\rm exc}}
\newcommand{\Wex}{{\rm Wex}}
\newcommand{\wex}{{\rm wex}}
\newcommand{\Fix}{{\rm Fix}}
\newcommand{\fix}{{\rm fix}}
\newcommand{\lrmax}{{\rm lrmax}}
\newcommand{\rlmax}{{\rm rlmax}}
\newcommand{\Rec}{{\rm Rec}}
\newcommand{\rec}{{\rm rec}}
\newcommand{\Arec}{{\rm Arec}}
\newcommand{\arec}{{\rm arec}}
\newcommand{\ERec}{{\rm ERec}}
\newcommand{\erec}{{\rm erec}}
\newcommand{\EArec}{{\rm EArec}}
\newcommand{\earec}{{\rm earec}}
\newcommand{\recarec}{{\rm recarec}}
\newcommand{\nonrec}{{\rm nonrec}}
\newcommand{\Cpeak}{{\rm Cpeak}}
\newcommand{\cpeak}{{\rm cpeak}}
\newcommand{\Cval}{{\rm Cval}}
\newcommand{\cval}{{\rm cval}}
\newcommand{\Cdasc}{{\rm Cdasc}}
\newcommand{\cdasc}{{\rm cdasc}}
\newcommand{\Cddes}{{\rm Cddes}}
\newcommand{\cddes}{{\rm cddes}}
\newcommand{\cdrise}{{\rm cdrise}}
\newcommand{\cdfall}{{\rm cdfall}}
\newcommand{\Peak}{{\rm Peak}}
\newcommand{\peak}{{\rm peak}}
\newcommand{\Val}{{\rm Val}}
\newcommand{\val}{{\rm val}}
\newcommand{\Dasc}{{\rm Dasc}}
\newcommand{\dasc}{{\rm dasc}}
\newcommand{\Ddes}{{\rm Ddes}}
\newcommand{\ddes}{{\rm ddes}}
\newcommand{\inv}{{\rm inv}}
\newcommand{\maj}{{\rm maj}}
\newcommand{\rs}{{\rm rs}}
\newcommand{\cross}{{\rm cr}}
\newcommand{\crosshat}{{\widehat{\rm cr}}}
\newcommand{\nest}{{\rm ne}}
\newcommand{\rodd}{{\rm rodd}}
\newcommand{\reven}{{\rm reven}}
\newcommand{\lodd}{{\rm lodd}}
\newcommand{\leven}{{\rm leven}}
\newcommand{\sg}{{\rm sg}}
\newcommand{\bl}{{\rm bl}}
\newcommand{\tran}{{\rm tr}}
\newcommand{\area}{{\rm area}}
\newcommand{\ret}{{\rm ret}}
\newcommand{\peaks}{{\rm peaks}}
\newcommand{\hl}{{\rm hl}}
\newcommand{\sll}{{\rm sl}}
\newcommand{\negg}{{\rm neg}}
\newcommand{\imp}{{\rm imp}}
\newcommand{\osg}{{\rm osg}}
\newcommand{\ons}{{\rm ons}}
\newcommand{\isg}{{\rm isg}}
\newcommand{\ins}{{\rm ins}}
\newcommand{\LL}{{\rm LL}}
\newcommand{\height}{{\rm ht}}
\newcommand{\as}{{\rm as}}

\newcommand{\ba}{{\bm{a}}}
\newcommand{\bahat}{{\widehat{\bm{a}}}}
\newcommand{\sfa}{{{\sf a}}}
\newcommand{\bb}{{\bm{b}}}
\newcommand{\bc}{{\bm{c}}}
\newcommand{\bchat}{{\widehat{\bm{c}}}}
\newcommand{\bd}{{\bm{d}}}
\newcommand{\bee}{{\bm{e}}}
\newcommand{\beh}{{\bm{eh}}}
\newcommand{\bff}{{\bm{f}}}
\newcommand{\bg}{{\bm{g}}}
\newcommand{\bh}{{\bm{h}}}
\newcommand{\bll}{{\bm{\ell}}}
\newcommand{\bp}{{\bm{p}}}
\newcommand{\br}{{\bm{r}}}
\newcommand{\bs}{{\bm{s}}}
\newcommand{\bu}{{\bm{u}}}
\newcommand{\bw}{{\bm{w}}}
\newcommand{\bx}{{\bm{x}}}
\newcommand{\by}{{\bm{y}}}
\newcommand{\bz}{{\bm{z}}}
\newcommand{\bA}{{\bm{A}}}
\newcommand{\bB}{{\bm{B}}}
\newcommand{\bC}{{\bm{C}}}
\newcommand{\bE}{{\bm{E}}}
\newcommand{\bF}{{\bm{F}}}
\newcommand{\bG}{{\bm{G}}}
\newcommand{\bH}{{\bm{H}}}
\newcommand{\bI}{{\bm{I}}}
\newcommand{\bJ}{{\bm{J}}}
\newcommand{\bL}{{\bm{L}}}
\newcommand{\bLhat}{{\widehat{\bm{L}}}}
\newcommand{\bM}{{\bm{M}}}
\newcommand{\bN}{{\bm{N}}}
\newcommand{\bP}{{\bm{P}}}
\newcommand{\bQ}{{\bm{Q}}}
\newcommand{\bR}{{\bm{R}}}
\newcommand{\bS}{{\bm{S}}}
\newcommand{\bT}{{\bm{T}}}
\newcommand{\bW}{{\bm{W}}}
\newcommand{\bX}{{\bm{X}}}
\newcommand{\bY}{{\bm{Y}}}
\newcommand{\bIB}{{\bm{B}^{\rm irr}}}
\newcommand{\bOB}{{\bm{B}^{\rm ord}}}
\newcommand{\bOS}{{\bm{OS}}}
\newcommand{\bERR}{{\bm{ERR}}}
\newcommand{\bSP}{{\bm{SP}}}
\newcommand{\bMV}{{\bm{MV}}}
\newcommand{\bBM}{{\bm{BM}}}
\newcommand{\balpha}{{\bm{\alpha}}}
\newcommand{\balphapre}{{\bm{\alpha}^{\rm pre}}}
\newcommand{\bbeta}{{\bm{\beta}}}
\newcommand{\bgamma}{{\bm{\gamma}}}
\newcommand{\bdelta}{{\bm{\delta}}}
\newcommand{\bkappa}{{\bm{\kappa}}}
\newcommand{\bmu}{{\bm{\mu}}}
\newcommand{\bomega}{{\bm{\omega}}}
\newcommand{\bsigma}{{\bm{\sigma}}}
\newcommand{\btau}{{\bm{\tau}}}
\newcommand{\bphi}{{\bm{\phi}}}
\newcommand{\bphihat}{{\skew{3}\widehat{\vphantom{t}\protect\smash{\bm{\phi}}}}}
   %% Moves the hat a bit down to make it look nicer
\newcommand{\bpsi}{{\bm{\psi}}}
\newcommand{\bxi}{{\bm{\xi}}}
\newcommand{\bzeta}{{\bm{\zeta}}}
\newcommand{\bone}{{\bm{1}}}
\newcommand{\bzero}{{\bm{0}}}

\newcommand{\Cbar}{{\overline{C}}}
\newcommand{\Dbar}{{\overline{D}}}
\newcommand{\dbar}{{\overline{d}}}
\def\Btilde{{\widetilde{B}}}
\def\Ctilde{{\widetilde{C}}}
\def\Ftilde{{\widetilde{F}}}
\def\Gtilde{{\widetilde{G}}}
\def\Htilde{{\widetilde{H}}}
\def\Lhat{{\widehat{L}}}
\def\Ltilde{{\widetilde{L}}}
\def\Ptilde{{\widetilde{P}}}
\def\ptilde{{\widetilde{p}}}
\def\Chat{{\widehat{C}}}
\def\ctilde{{\widetilde{c}}}
\def\zbar{{\overline{Z}}}
\def\pitilde{{\widetilde{\pi}}}
\def\omegahat{{\widehat{\omega}}}

\newcommand{\sech}{{\rm sech}}

%
% Jacobian and Dixonian elliptic functions
%
\newcommand{\sn}{{\rm sn}}
\newcommand{\cn}{{\rm cn}}
\newcommand{\dn}{{\rm dn}}
\newcommand{\sm}{{\rm sm}}
\newcommand{\cm}{{\rm cm}}

%
% Commands for hypergeometric series
%
\newcommand{\zfz}{ {{}_0 \! F_0} }
\newcommand{\zfo}{ {{}_0  F_1} }
\newcommand{\ofz}{ {{}_1 \! F_0} }
\newcommand{\ofo}{ {{}_1 \! F_1} }
\newcommand{\oft}{ {{}_1 \! F_2} }
%\newcommand{\tfo}{ {{}_2 \! F_1} }

%
% Hypergeometric functions, using "tensor" package
%
\newcommand{\FHyper}[2]{ {\tensor[_{#1 \!}]{F}{_{#2}}\!} }
\newcommand{\FHYPER}[5]{ {\FHyper{#1}{#2} \!\biggl(
   \!\!\begin{array}{c} #3 \\[1mm] #4 \end{array}\! \bigg|\, #5 \! \biggr)} }
\newcommand{\tfo}{ {\FHyper{2}{1}} }
\newcommand{\tfz}{ {\FHyper{2}{0}} }
\newcommand{\threefz}{ {\FHyper{3}{0}} }
\newcommand{\FHYPERbottomzero}[3]{ {\FHyper{#1}{0} \hspace*{-0mm}\biggl(
   \!\!\begin{array}{c} #2 \\[1mm] \hbox{---} \end{array}\! \bigg|\, #3 \! \biggr)} }
\newcommand{\FHYPERtopzero}[3]{ {\FHyper{0}{#1} \hspace*{-0mm}\biggl(
   \!\!\begin{array}{c} \hbox{---} \\[1mm] #2 \end{array}\! \bigg|\, #3 \! \biggr)} }

\newcommand{\phiHyper}[2]{ {\tensor[_{#1}]{\phi}{_{#2}}} }
\newcommand{\psiHyper}[2]{ {\tensor[_{#1}]{\psi}{_{#2}}} }
\newcommand{\PhiHyper}[2]{ {\tensor[_{#1}]{\Phi}{_{#2}}} }
\newcommand{\PsiHyper}[2]{ {\tensor[_{#1}]{\Psi}{_{#2}}} }
\newcommand{\phiHYPER}[6]{ {\phiHyper{#1}{#2} \!\left(
   \!\!\begin{array}{c} #3 \\ #4 \end{array}\! ;\, #5, \, #6 \! \right)\!} }
\newcommand{\psiHYPER}[6]{ {\psiHyper{#1}{#2} \!\left(
   \!\!\begin{array}{c} #3 \\ #4 \end{array}\! ;\, #5, \, #6 \! \right)} }
\newcommand{\PhiHYPER}[5]{ {\PhiHyper{#1}{#2} \!\left(
   \!\!\begin{array}{c} #3 \\ #4 \end{array}\! ;\, #5 \! \right)\!} }
\newcommand{\PsiHYPER}[5]{ {\PsiHyper{#1}{#2} \!\left(
   \!\!\begin{array}{c} #3 \\ #4 \end{array}\! ;\, #5 \! \right)\!} }
\newcommand{\zerophizero}{ {\phiHyper{0}{0}} }
\newcommand{\ophizero}{ {\phiHyper{1}{0}} }
\newcommand{\zphio}{ {\phiHyper{0}{1}} }
\newcommand{\ophio}{ {\phiHyper{1}{1}} }
\newcommand{\tphio}{ {\phiHyper{2}{1}} }
\newcommand{\tphiz}{ {\phiHyper{2}{0}} }
\newcommand{\tPhio}{ {\PhiHyper{2}{1}} }
\newcommand{\opsio}{ {\psiHyper{1}{1}} }

%
% Variants of \binom  (defined using the AMS "genfrac" command)
%
\newcommand{\stirlingsubset}[2]{\genfrac{\{}{\}}{0pt}{}{#1}{#2}}
\newcommand{\stirlingcycle}[2]{\genfrac{[}{]}{0pt}{}{#1}{#2}}
\newcommand{\assocstirlingsubset}[3]{{\genfrac{\{}{\}}{0pt}{}{#1}{#2}}_{\! \ge #3}}
\newcommand{\genstirlingsubset}[4]{{\genfrac{\{}{\}}{0pt}{}{#1}{#2}}_{\! #3,#4}}
\newcommand{\irredstirlingsubset}[2]{{\genfrac{\{}{\}}{0pt}{}{#1}{#2}}^{\!\rm irr}}
\newcommand{\euler}[2]{\genfrac{\langle}{\rangle}{0pt}{}{#1}{#2}}
\newcommand{\eulergen}[3]{{\genfrac{\langle}{\rangle}{0pt}{}{#1}{#2}}_{\! #3}}
\newcommand{\eulersecond}[2]{\left\langle\!\! \euler{#1}{#2} \!\!\right\rangle}
\newcommand{\eulersecondgen}[3]{{\left\langle\!\! \euler{#1}{#2} \!\!\right\rangle}_{\! #3}}
\newcommand{\binomvert}[2]{\genfrac{\vert}{\vert}{0pt}{}{#1}{#2}}
\newcommand{\binomsquare}[2]{\genfrac{[}{]}{0pt}{}{#1}{#2}}
\newcommand{\doublebinom}[2]{\left(\!\! \binom{#1}{#2} \!\!\right)}
\newcommand{\lahnum}[2]{\genfrac{\lfloor}{\rfloor}{0pt}{}{#1}{#2}}
\newcommand{\rlahnum}[3]{{\genfrac{\lfloor}{\rfloor}{0pt}{}{#1}{#2}}_{\! #3}}

%
% Definitions for Laguerre and Lah polynomials
%
\newcommand{\Lna}{{L_n^{(\alpha)}}}
\newcommand{\bfLna}{{\bfL_n^{(\alpha)}}}
\newcommand{\bfLnamult}{{\bfL^{(\balpha)}_{\bfn}}}
\newcommand{\scrlna}{{\scrl_n^{(\alpha)}}}
\newcommand{\scrlmultna}{{\scrl_\bfn^{(\balpha)}}}
\newcommand{\scrlbarna}{{\scrlbar_n^{(\alpha)}}}
\newcommand{\Lah}{{\textrm{Lah}}}
\newcommand{\LD}{{\mathbf{LD}}}

% Array for subscripts

\newenvironment{sarray}{
             \textfont0=\scriptfont0
             \scriptfont0=\scriptscriptfont0
             \textfont1=\scriptfont1
             \scriptfont1=\scriptscriptfont1
             \textfont2=\scriptfont2
             \scriptfont2=\scriptscriptfont2
             \textfont3=\scriptfont3
             \scriptfont3=\scriptscriptfont3
           \renewcommand{\arraystretch}{0.7}
           \begin{array}{l}}{\end{array}}

\newenvironment{scarray}{
             \textfont0=\scriptfont0
             \scriptfont0=\scriptscriptfont0
             \textfont1=\scriptfont1
             \scriptfont1=\scriptscriptfont1
             \textfont2=\scriptfont2
             \scriptfont2=\scriptscriptfont2
             \textfont3=\scriptfont3
             \scriptfont3=\scriptscriptfont3
           \renewcommand{\arraystretch}{0.7}
           \begin{array}{c}}{\end{array}}

% Circled math symbols:
% From http://latex-community.org/forum/viewtopic.php?f=44&t=22367

%\usepackage{tikz}
\newcommand*\circled[1]{\tikz[baseline=(char.base)]{
  \node[shape=circle,draw,inner sep=1pt] (char) {#1};}}
\newcommand{\ostar}{{\circledast}}
\newcommand{\ostarN}{{\,\circledast_{\vphantom{\dot{N}}N}\,}}
\newcommand{\ostarPsi}{{\,\circledast_{\vphantom{\dot{\Psi}}\Psi}\,}}
\newcommand{\starN}{{\,\ast_{\vphantom{\dot{N}}N}\,}}
\newcommand{\starpsi}{{\,\ast_{\vphantom{\dot{\bpsi}}\!\bpsi}\,}}
\newcommand{\starone}{{\,\ast_{\vphantom{\dot{1}}1}\,}}
\newcommand{\startwo}{{\,\ast_{\vphantom{\dot{2}}2}\,}}
\newcommand{\starinfty}{{\,\ast_{\vphantom{\dot{\infty}}\infty}\,}}
\newcommand{\starT}{{\,\ast_{\vphantom{\dot{T}}T}\,}}

%% For scaling equations (uses "graphicx" package):  see
%% http://tex.stackexchange.com/questions/60453/reducing-font-size-in-equation
\newcommand*{\Scale}[2][4]{\scalebox{#1}{$#2$}}

\newcommand*{\Scaletext}[2][4]{\scalebox{#1}{#2}} %% THIS DOESN'T SEEM TO WORK

%% \clearpage 

%% \tableofcontents

\clearpage

\section{Introduction}

The monic Laguerre polynomials $\bfLna(x) = (-1)^n n! \, \Lna(x)$
can be defined as \cite{Rainville_60,Szego_75,Andrews_99,Ismail_05}
\begin{subeqnarray}
  \bfLna(x)
  & = &
  (-1)^n \, (\alpha+1)^{\overline{n}} \; \FHYPER{1}{1}{-n}{\alpha+1}{x}
      \slabel{def.bfLna.a}  \\[2mm]
  & = &
  \sum_{k=0}^n (-1)^{n-k} \, \binom{n}{k} \, (\alpha+1+k)^{\overline{n-k}} \, x^k
      \slabel{def.bfLna.b}
 \label{def.bfLna}
\end{subeqnarray}
where $r^{\overline{n}} \eqdef r(r+1) \cdots (r+n-1)$;
note that they are polynomials (with integer coefficients)
jointly in $x$ and $\alpha$.
The monic Laguerre polynomials have the exponential generating function
\be
  \sum_{n=0}^\infty \bfLna(x) \: {t^n \over n!}
  \;=\;
  (1+t)^{-(\alpha+1)} \, e^{xt/(1+t)}
  \;.
 \label{eq.genfn.1}
\ee
For $\alpha > -1$ they are orthogonal with respect to the measure
$x^\alpha e^{-x} \, dx$ on $(0,\infty)$.
Using Kummer's first transformation
for the confluent hypergeometric function $\ofo$
\cite[eq.~(1.4.11)]{Ismail_05},
eq.~\reff{def.bfLna.a} can also be rewritten as
\be
  \bfLna(x)
  \;=\;
  (-1)^n \, (\alpha+1)^{\overline{n}} \, e^x \;
       \FHYPER{1}{1}{\alpha+1+n}{\alpha+1}{-x}
  \;.
 \label{def.bfLna.bis}
\ee

Now fix an integer $r \ge 1$.
The multiple Laguerre polynomials of the first kind of type~II
\cite[section~23.4.1]{Ismail_05},
denoted $\bfLnamult(x)$
%% denoted $\bfL^{(\balpha)}_{\bfn}(x)$
where $\balpha = (\alpha_1,\ldots,\alpha_r)$ and $\bfn = (n_1,\ldots,n_r)$,
can be defined by a straightforward generalization of \reff{def.bfLna.bis}:
\be
  \bfLnamult(x)
  \;=\;
  (-1)^{|\bfn|} \,
       \biggl( \prod\limits_{i=1}^r (\alpha_i+1)^{\overline{n_i}} \biggr)
       \, e^x \;
       \FHYPER{r}{r}{\alpha_1 + 1 + n_1 ,\,\ldots,\, \alpha_r + 1 + n_r}
                    {\alpha_1 + 1 ,\,\ldots,\, \alpha_r + 1}{-x}
 \label{def.bfLnamult.bis}
\ee
where $|\bfn| \eqdef n_1 + \ldots + n_r$.
It follows from known properties of the hypergeometric function $\FHyper{r}{r}$
that the right-hand side of \reff{def.bfLnamult.bis}
is an entire function of $x$
that behaves asymptotically at infinity like $x^{|\bfn|}$;
%% $x^{n_1 + \ldots + n_r}$;
therefore it is a (monic) polynomial in $x$,
of degree $|\bfn|$.\footnote{
   This reasoning goes back at least to Hille \cite[p.~52]{Hille_29}.
   The needed asymptotic expansion of $\FHyper{r}{r}$
   can be found in \cite[section~5.11.3]{Luke_69}
   or \cite{Volkmer_14}.
}
In fact, we have the explicit expression,
which generalizes \reff{def.bfLna.b}:\footnote{
   This formula follows from \reff{def.bfLnamult.bis}
   by application of Karlsson's \cite{Karlsson_71} identity
   for hypergeometric functions where the numerator and denominator parameters
   differ by integers, combined with
   %% $\FHYPER{0}{0}{\hbox{---}}{\hbox{---}}{-\!x} = e^{-x}$
   %% $\FHyper{0}{0}(\,\hbox{---}\, | -\!x) = e^{-x}$
   $\FHyper{0}{0}(\!\begin{array}{cc}\hbox{---}\\[-2mm] \hbox{---} \end{array}\! | -\!x) = e^{-x}$
    at the final stage.
   See also Srivastava \cite{Srivastava_73} for a very simple proof
   of Karlsson's identity;
   and see \cite{Chakrabarty_74,Panda_76} for some interesting generalizations.
}
% \begin{eqnarray}
%    & &
%    \!\!\!\!\!\!
%    \bfLnamult(x)
%    \;=\;
%    \sum_{k_1 = 0}^{n_1} \cdots \sum_{k_r = 0}^{n_r}
%       (-1)^{|\bfn| - |\bfk|} \:
%       \binom{n_1}{k_1} \,\cdots\, \binom{n_r}{k_r} \:
%       \times\:
%              \nonumber \\[1mm]
%    & &
%    (\alpha_1 + 1 + k_1 + \ldots + k_r)^{\overline{n_1 - k_1}}
%    (\alpha_2 + 1 + k_2 + \ldots + k_r)^{\overline{n_2 - k_2}}
%    \,\cdots\,
%    (\alpha_2 + 1 + k_r)^{\overline{n_r - k_r}}
%    \:
%    x^{|\bfk|}
%    \:.
%    \qquad
%  \label{def.bfLnamult.explicit}
% \end{eqnarray}
\be
   \bfLnamult(x)
   \;=\;
   \sum_{k_1 = 0}^{n_1} \cdots \sum_{k_r = 0}^{n_r}
      (-1)^{|\bfn| - |\bfk|} \:
      \Biggl( \prod_{i=1}^r \binom{n_i}{k_i} \;
%%          (\alpha_i + 1 + k_i + k_{i+1} + \ldots + k_r)^{\overline{n_i - k_i}}
            (\alpha_i + 1 + k_1 + \ldots + k_i)^{\overline{n_i - k_i}}
      \Biggr)
      \: x^{|\bfk|}
   \;.
 \label{def.bfLnamult.explicit}
\ee
When $\alpha_1,\ldots,\alpha_r > -1$
with $\alpha_i - \alpha_j \notin \Z$ for all pairs $i \neq j$,
these polynomials are multiple orthogonal \cite[Chapter~23]{Ismail_05}
with respect to the collection of measures
$x^{\alpha_i} e^{-x} \, dx$ on $(0,\infty)$ with $1 \le i \le r$.
Finally, the multiple Laguerre polynomials have the
multivariate exponential generating function
\cite{Lee_07}
\be
   \sum_{n_1 = 0}^\infty \cdots \sum_{n_r = 0}^\infty
   \bfLnamult(x) \:
   {t_1^{n_1} \over n_1!} \,\cdots\, {t_r^{n_r} \over n_r!}
   \;=\;
   \biggl( \prod\limits_{i=1}^r (1+t_i)^{-(\alpha_i+1)} \biggr)
   \exp\biggl[ x \Bigl( 1 - \prod\limits_{i=1}^r {1 \over 1 + t_i} \Bigr)
       \biggr]
   \;.
 \label{def.bfLnamult.EGF}
\ee

\bigskip

{\bf Remark/Question.}  The multiple Laguerre polynomial $\bfLnamult(x)$
is invariant under joint permutations of $\bfn$ and $\balpha$:
this is manifest in \reff{def.bfLnamult.bis} and \reff{def.bfLnamult.EGF},
but is far from obvious in the explicit formula \reff{def.bfLnamult.explicit}.
Is there some easy way of deriving this symmetry from
\reff{def.bfLnamult.explicit}?
And is there an alternate explicit formula in which this symmetry is manifest?
\myendremark

%% {\bf CHECK ALSO MULTIPLE LAGUERRE POLYNOMIALS OF THE \emph{SECOND} KIND!!!}

The purpose of the present paper is twofold:
(a) to give a combinatorial interpretation of the multiple Laguerre polynomials
\reff{def.bfLnamult.bis}/\reff{def.bfLnamult.explicit},
generalizing the digraph model found by
Foata and Strehl \cite{Foata_84} for the ordinary Laguerre polynomials;
and (b) to give an explicit integral representation for these polynomials,
showing that they form a multidimensional Stieltjes moment sequence
whenever $x \le 0$.

\section{Combinatorial model}

Three decades ago, Foata and Strehl \cite{Foata_84}
introduced a beautiful combinatorial interpretation
of the Laguerre polynomials.
Let us define a \textbfit{Laguerre digraph}
to be a digraph in which each vertex has out-degree 0 or 1
and in-degree 0 or 1.
It follows that each weakly connected component
% \footnote{
%    By a ``weakly connected component'' of a digraph,
%    we mean a connected component of the underlying undirected graph.
%% THIS SHOULD BE "CORRESPONDING TO" a connected component of the
%%    underlying undirected graph.
% }
of a Laguerre digraph
is either a directed path of some length $\ell \ge 0$
(where a path of length 0 is an isolated vertex)
or else a directed cycle of some length $\ell \ge 1$
(where a cycle of length 1 is a loop).
Let us write $\LD_n$ for the set of Laguerre digraphs
on the vertex set $[n] \eqdef \{1,\ldots,n\}$;
and for a Laguerre digraph $G$,
let us write $\cyc(G)$ [resp.\ $\pa(G)$]
for the number of cycles (resp.\ paths) in $G$.
Foata and Strehl \cite{Foata_84} then showed that
the monic unsigned Laguerre polynomials
\be
   \scrlna(x)
   \;\eqdef\;
   n! \, \Lna(-x)
   \;=\;
   (-1)^n \, \bfLna(-x)
 \label{def.scrlna}
\ee
have the combinatorial representation
\be
   \scrlna(x)
   \;=\;
   \sum_{G \in \LD_n}  x^{\pa(G)} \, (\alpha+1)^{\cyc(G)}
   \;.
 \label{eq.foata.0}
\ee
Indeed, the proof of \reff{eq.foata.0} is an easy argument
using the exponential formula \cite[chapter~5]{Stanley_99},
or equivalently, the theory of species \cite{Bergeron_98}:
the number of directed paths on $n \ge 1$ vertices is $n!$,
so with a weight $x$ per path they have exponential generating function
$xt/(1-t)$.
The number of directed cycles on $n \ge 1$ vertices is $(n-1)!$,
so with a weight $\alpha+1$ per cycle
they have exponential generating function $- (\alpha+1) \log(1-t)$.
A Laguerre digraph is a disjoint union of paths and cycles,
so by the exponential formula it has exponential generating function
\be
   \exp\Bigl[ {xt \over 1-t} \,-\, (\alpha+1) \log(1-t) \Bigr]
   \;=\;
   (1-t)^{-(\alpha+1)} \, e^{xt/(1-t)}
   \;,
\ee
which coincides with \reff{eq.genfn.1} after $x \to -x$ and $t \to -t$.
Foata and Strehl \cite{Foata_84} also gave a direct combinatorial proof
of \reff{eq.foata.0} based on the definition \reff{def.bfLna};
this requires a bit more work \cite[Lemma~2.1]{Foata_84}.

Our first result is a combinatorial interpretation
of the multiple Laguerre polynomials
that extends the Foata--Strehl interpretation to $r > 1$.
For $\bfn = (n_1,\ldots,n_r) \in \N^r$,
we define a digraph $G_\bfn = (V_\bfn, \vec{E}_\bfn)$ with vertex set
\be
   V_\bfn
   \;=\;  
   \{ (i,j) \colon\: 1 \le i \le r \hbox{ and } 1 \le j \le n_i \}
\ee
and edge set
\be
   \vec{E}_\bfn
   \;=\;
   \bigl\{ \overrightarrow{(i,j) \, (i',j')} \colon\:  i \le i' \bigr\}
   \;.
\ee
The vertex set is thus the disjoint union of ``layers'' $V_i \simeq [n_i]$
for $1 \le i \le r$;
the edge set consists of all possible directed edges (including loops)
within each layer $V_i$, together with all possible edges from a layer $V_i$
to a layer $V_{i'}$ with $i' > i$.
We then write $\LD_\bfn$ for the set of Laguerre digraphs
that are spanning subdigraphs of $G_\bfn$,
i.e.\ Laguerre digraphs of the form
$(V_\bfn,A)$ with $A \subseteq \vec{E}_\bfn$.
Note that in a Laguerre digraph $G \in \LD_\bfn$,
every cycle must lie in a single layer $V_i$;
we denote by $\cyc_i(G)$ the number of cycles in layer $V_i$.
We then have:

\begin{theorem}
   \label{thm1.1}
The monic unsigned multiple Laguerre polynomials
\be
   \scrlmultna(x)
   \;\eqdef\;
   (-1)^{|\bfn|} \, \bfLnamult(-x)
 \label{def.scrlmultna}
\ee
have the combinatorial representation
\be
   \scrlmultna(x)
   \;=\;
   \sum_{G \in \LD_\bfn}
       x^{\pa(G)} \, \prod\limits_{i=1}^r (\alpha_i+1)^{\cyc_i(G)}
   \;.
 \label{eq.thm1.1}
\ee
\end{theorem}

The proof of this result is a simple generalization
of the argument just given for the Foata--Strehl formula \reff{eq.foata.0}:
%% indeed, it is so simple that we give it now:

\proofof{Theorem~\ref{thm1.1}}
Denote the right-hand side of \reff{eq.thm1.1} by
$\widehat{\scrl}_\bfn^{(\balpha)}(x)$,
and consider its multivariate exponential generating function
\be
   F(t_1,\ldots,t_r)
   \;\eqdef\;
   \sum_{n_1 = 0}^\infty \cdots \sum_{n_r = 0}^\infty
   \widehat{\scrl}_\bfn^{(\balpha)}(x)  \:
   {t_1^{n_1} \over n_1!} \,\cdots\, {t_r^{n_r} \over n_r!}
   \;.
 \label{def.scrlnamult.EGF}
\ee
We again argue using the exponential formula.
The multivariate exponential generating function
for a single directed cycle in layer $V_i$ is,
as before, $- (\alpha_i+1) \log(1-t_i)$.
Let~us now look at paths.
Every path $P$ in the digraph $G_\bfn$ is of the following form:
In~each layer $V_i$ choose a directed path $P_i$;
the $P_i$ are allowed to be empty, provided that they are not {\em all}\/ empty.
Let $i_1 < i_2 < \ldots < i_k$ be the indices with $P_i$ nonempty,
and construct the path $P$ obtained from the union of the $P_i$
by adjoining the edge linking
the final vertex of $P_{i_1}$ to the initial vertex of $P_{i_2}$,
the edge linking
the final vertex of $P_{i_2}$ to the initial vertex of $P_{i_3}$, etc.
With a weight $x$ per path, the multivariate exponential generating function
for a single such path is
\be
   x \, \Bigl( \prod\limits_{i=1}^r {1 \over 1 - t_i} \:-\:  1 \Bigr)
   \;.
\ee
Therefore, by the exponential formula we have
\be
   F(t_1,\ldots,t_r)
   \;=\;
   \exp\biggl[ - \sum_{i=1}^r (\alpha_i+1) \log(1-t_i)
               \:+\:
   x \, \Bigl( \prod\limits_{i=1}^r {1 \over 1 - t_i} \:-\:  1 \Bigr)
       \biggr]
   \;,
\ee
which coincides with \reff{def.bfLnamult.EGF}
after $x \to -x$ and $t_i \to -t_i$.
\qed

{\bf Remarks.}
1.  We leave it as an open problem to devise a direct combinatorial proof
of \reff{eq.thm1.1} based on the explicit formula \reff{def.bfLnamult.explicit}.

2.  The combinatorial representation \reff{eq.thm1.1},
unlike the explicit formula \reff{def.bfLnamult.explicit},
manifestly exhibits the invariance of $\scrlmultna(x)$
under joint permutations of $\bfn$ and $\balpha$,
since there is a weight-preserving bijection between
the digraphs contributing to the right-hand side of \reff{eq.thm1.1}
for the original and permuted cases.
I thank an anonymous referee for pointing this out.

3.  For the case $r=2$, a slightly different combinatorial interpretation
of the multiple Laguerre polynomials was found by
Drake \cite[Theorem~3.5.2]{Drake_06}.
But also this representation fails to manifestly exhibit
the permutation symmetry.
\myendremark

\section{Stieltjes moment representation}

For the ordinary Laguerre polynomials ($r=1$),
a well-known integral representation \cite[Theorem~5.4]{Szego_75}
asserts that
\be
   \scrlna(x)
   \;=\;
   n! \, L_n^{(\alpha)}(-x)
   \;=\;
   e^{-x} x^{-\alpha/2}
   \int\limits_0^\infty y^n \: e^{-y} \, y^{\alpha/2} \,
         I_\alpha(2 \sqrt{xy}) \: dy
   \qquad\hbox{for } \alpha > -1
   \;,
 \label{eq.laguerre.stieltjes}
\ee
where $I_\alpha$ is a modified Bessel function
of the first kind \cite[p.~77]{Watson_44}:
\begin{subeqnarray}
   I_\alpha(z)
   & = &
   \sum_{k=0}^\infty {(z/2)^{\alpha+2k} \over k! \, \Gamma(\alpha+k+1)}
        \\[2mm]
   & = &
   \displaystyle {1 \over \Gamma(\alpha+1)} \: (z/2)^\alpha \: 
      \FHYPERtopzero{1}{\alpha+1}{z^2/4}
   \;.
 \label{def.BesselI}
\end{subeqnarray}
Since $I_\alpha$ is nonnegative on $[0,\infty)$,
it follows from \reff{eq.laguerre.stieltjes}
that the sequence $(\scrlna(x))_{n \ge 0}$
is a Stieltjes moment sequence whenever $\alpha \ge -1$ and $x \ge 0$:
that is,
\be
   \scrlna(x)
   \;=\;
   \int\limits_0^\infty y^n \: d\mu_{\alpha,x}(y)
\ee
where
\be
   d\mu_{\alpha,x}(y)
   \;=\;
   \begin{cases}
     e^{-x} \; \FHYPERtopzero{1}{\alpha+1}{xy}
            \: \displaystyle {1 \over \Gamma(\alpha+1)}
            \: y^\alpha \, e^{-y} \, dy
         & \textrm{for $\alpha > -1$}  \\[6mm]
     x \, e^{-(x+y)} \; \FHYPERtopzero{1}{2}{xy} \, dy
         & \textrm{for $\alpha = -1$}
   \end{cases}
\ee
is a positive measure on $[0,\infty)$.\footnote{
   All this was observed a half-century ago by Karlin
   \cite[p.~62]{Karlin_68b} \cite[pp.~440--441]{Karlin_68}.
}

We now give an integral representation for the
multiple Laguerre polynomials that generalizes \reff{eq.laguerre.stieltjes}
to $r > 1$:

\begin{theorem}
   \label{thm1.2}
Let $\alpha_1,\ldots,\alpha_r \ge -1$ and $x \ge 0$.
Then the multisequence $(\scrlmultna(x))_{\bfn \in \N^r}$
of monic unsigned multiple Laguerre polynomials
is a multidimensional Stieltjes moment sequence:
that is, there exists a positive measure $\mu_{\balpha,x}$ on $[0,\infty)^r$
such that
\be
   \scrlmultna(x)
   \;=\;
   \int\limits_{[0,\infty)^r}  \!\!\! \bfy^\bfn \: d\mu_{\balpha,x}(\bfy)
\ee
for all $\bfn \in \N^r$,
where $\bfy^\bfn \eqdef \prod\limits_{i=1}^r y_i^{n_i}$.
In fact, for $\alpha_1,\ldots,\alpha_r > -1$ we have the explicit formula
\be
   d\mu_{\balpha,x}(\bfy)
   \;=\;
   e^{-x} \;
     \FHYPERtopzero{r}{\alpha_1+1,\,\ldots,\, \alpha_r+1}{xy_1 \cdots y_r} \:
     \prod\limits_{i=1}^r  {1 \over \Gamma(\alpha_i+1)} \:
                              y_i^{\alpha_i} \, e^{-y_i} \, dy_i
   \;.
   \quad
 \label{eq.thm1.2}
\ee
\end{theorem}

\proof
We begin from the exponential generating function \reff{def.bfLnamult.EGF}
with $x \to -x$:
\be
   e^{-x} \:
   \biggl( \prod\limits_{i=1}^r (1+t_i)^{-(\alpha_i+1)} \biggr)
      \exp\biggl[ x \prod\limits_{i=1}^r {1 \over 1 + t_i} \biggr]
   \;=\;
   e^{-x} \:
   \sum_{n=0}^\infty {x^n \over n!}
                     \prod\limits_{i=1}^r (1+t_i)^{-(\alpha_i+1+n)}
   \:.
   \quad
\ee
We now assume that $\alpha_1,\ldots,\alpha_r > -1$
and insert the integral representation
\be
   (1+t_i)^{-(\alpha_i+1+n)}
   \;=\;
   {1 \over \Gamma(\alpha_i+1+n)} \:
     \int\limits_0^\infty e^{-t_i y_i} \, y_i^{n+\alpha_i} \, e^{-y_i} \, dy_i
   \;.
\ee
It follows that
\be
   e^{-x} \:
   \biggl( \prod\limits_{i=1}^r (1+t_i)^{-(\alpha_i+1)} \biggr)
      \exp\biggl[ x \prod\limits_{i=1}^r {1 \over 1 + t_i} \biggr]
   \;=\;
   \int\limits_{[0,\infty)^r} e^{-\bft\cdot\bfy} \: d\mu_{\balpha,x}(\bfy)
\ee
where
\begin{subeqnarray}
   d\mu_{\balpha,x}(\bfy)
   & = &
   e^{-x} \:
   \sum\limits_{n=0}^\infty {x^n \over n!}
       \prod\limits_{i=1}^r
          {y_i^{n+\alpha_i} \, e^{-y_i} \over \Gamma(\alpha_i+1+n)} \: dy_i
     \\[2mm]
   & = &
   e^{-x} \;
     \FHYPERtopzero{r}{\alpha_1+1,\,\ldots,\, \alpha_r+1}{xy_1 \cdots y_r} \:
     \prod\limits_{i=1}^r  {1 \over \Gamma(\alpha_i+1)} \:
                              y_i^{\alpha_i} \, e^{-y_i} \, dy_i
   \;.
      \nonumber \\[-3mm]
\end{subeqnarray}
Extracting the coefficient of $\bft^\bfn/\bfn!$, we conclude that
\be
   \scrlmultna(x)
   \;=\;
   \int\limits_{[0,\infty)^r}  \!\!\! \bfy^\bfn \: d\mu_{\balpha,x}(\bfy)
   \;.
\ee
This shows that $(\scrlmultna(x))_{\bfn \in \N^r}$
is a multidimensional Stieltjes moment sequence
whenever $\alpha_1,\ldots,\alpha_r > -1$;
and it holds also for $\alpha_1,\ldots,\alpha_r \ge -1$
since the set of multidimensional Stieltjes moment sequences
is closed under pointwise limits.
\qed

%% \noindent
%% We will prove Theorem~\ref{thm1.2} in Section~\ref{sec.proof.stieltjes};
%% it is an easy application of the theory of
%% multidimensional complete monotonicity.

In particular, Theorem~\ref{thm1.2} implies:

\begin{corollary}
   \label{cor.thm1.2}
Let $\alpha_1,\ldots,\alpha_r \ge -1$ and $x \ge 0$,
and fix a multi-index $\bfk \in \N^r$.
Then the sequence $(\scrl_{n\bfk}^{(\balpha)}(x))_{n \ge 0}$
is a Stieltjes moment sequence:
that is, there exists a positive measure $\mu_{\balpha,x,\bfk}$ on $[0,\infty)$
such that
\be
   \scrl_{n\bfk}^{(\balpha)}(x)
   \;=\;
   \int\limits_{[0,\infty)}  \!\! y^n \: d\mu_{\balpha,x,\bfk}(y)
\ee
for all $n \ge 0$.
\end{corollary}

Corollary~\ref{cor.thm1.2} can be restated in the language of total positivity.
Recall that a finite or infinite matrix of real numbers is called
{\em totally positive}\/ (TP) if all its minors are nonnegative,
and {\em totally positive of order~$r$} (TP${}_r$)
if all its minors of size $\le r$ are nonnegative.
Background information on totally positive matrices can be found
in \cite{Karlin_68,Gantmacher_02,Pinkus_10,Fallat_11};
they have application to many fields of pure and applied mathematics.
In particular, it is known
\cite[Th\'eor\`eme~9]{Gantmakher_37} \cite[section~4.6]{Pinkus_10}
that an infinite Hankel matrix $(a_{i+j})_{i,j \ge 0}$
of real numbers is totally positive if and only if the underlying sequence
$(a_n)_{n \ge 0}$ is a Stieltjes moment sequence.
%% i.e.\ the moments of a positive measure on $[0,\infty)$.
So Corollary~\ref{cor.thm1.2} asserts that, for every $\bfk \in \N^r$,
every minor of the infinite Hankel matrix
$(\scrl_{(i+j)\bfk}^{(\balpha)}(x))_{i,j \ge 0}$
is a polynomial in $x$ and $\alpha_1,\ldots,\alpha_r$
that is nonnegative whenever $\alpha_1,\ldots,\alpha_r \ge -1$ and $x \ge 0$.

But much more appears to be true:  namely, it seems that we have
{\em coefficientwise}\/ Hankel-total positivity
\cite{Sokal_flajolet,Sokal_OPSFA,Sokal_totalpos}
in the variables $x$ and $\beta_i \eqdef \alpha_i + 1$:

\begin{conjecture}[Coefficientwise Hankel-total positivity of the multiple Laguerre polynomials]
   \label{conj.totalpos}
For each multi-index $\bfk \in \N^r$,
the sequence $\big( \scrl_{n\bfk}^{(\bbeta-\bone)}(x) \big)_{n \ge 0}$
is coefficientwise Hankel-totally positive
in the variables $x$ and $\bbeta = (\beta_1,\ldots,\beta_r)$:
that is, every minor of the infinite Hankel matrix
$\big( \scrl_{(i+j)\bfk}^{(\bbeta-\bone)}(x) \big)_{i,j \ge 0}$
is a polynomial in $x$ and~$\bbeta$ with nonnegative coefficients.
\end{conjecture}

By symbolic computation using {\sc Mathematica},
I have verified this conjecture for the following cases:
\begin{itemize}
   \item $r=1$ and $\bfk = (1)$ up to the $11 \times 11$ Hankel matrix;
\\[-8mm]
   \item $r=2$ and $\bfk = (1,1)$ up to the $9 \times 9$ Hankel matrix;
\\[-8mm]
   \item $r=2$ and $\bfk = (2,1)$ up to the $8 \times 8$ Hankel matrix;
\\[-8mm]
   \item $r=2$ and $\bfk = (3,1)$ up to the $8 \times 8$ Hankel matrix;
\\[-8mm]
   \item $r=2$ and $\bfk = (3,2)$ up to the $8 \times 8$ Hankel matrix;
\\[-8mm]
   \item $r=3$ and $\bfk = (1,1,1)$ up to the $7 \times 7$ Hankel matrix;
\\[-8mm]
   \item $r=3$ and $\bfk = (2,1,1)$ up to the $6 \times 6$ Hankel matrix;
\\[-8mm]
   \item $r=3$ and $\bfk = (2,2,1)$ up to the $6 \times 6$ Hankel matrix;
\\[-8mm]
   \item $r=4$ and $\bfk = (1,1,1,1)$ up to the $6 \times 6$ Hankel matrix;
\\[-8mm]
   \item $r=4$ and $\bfk = (2,1,1,1)$ up to the $5 \times 5$ Hankel matrix;
\\[-8mm]
   \item $r=5$ and $\bfk = (1,1,1,1,1)$ up to the $4 \times 4$ Hankel matrix.
\end{itemize}
%% totalpos_hankel_laguerre.out (and out2, out3) on sokal.physics.nyu.edu
%% totalpos_hankel_higherlaguerre_p=2.out on sokal.physics.nyu.edu
%% totalpos_hankel_multiplelaguerre_21.out on sokal.physics.nyu.edu
%% totalpos_hankel_multiplelaguerre_31.out on sokal.physics.nyu.edu
%% totalpos_hankel_multiplelaguerre_32.out on sokal.physics.nyu.edu
%% totalpos_hankel_higherlaguerre_p=3.out on sokal.physics.nyu.edu
%% totalpos_hankel_multiplelaguerre_211.out on sokal.physics.nyu.edu
%% totalpos_hankel_multiplelaguerre_221.out on sokal.physics.nyu.edu
%% totalpos_hankel_higherlaguerre_p=4.out on sokal.physics.nyu.edu
%% totalpos_hankel_multiplelaguerre_2111.out on sokal.physics.nyu.edu
%% totalpos_hankel_higherlaguerre_p=5.out on sokal.physics.nyu.edu
%
For the case of ordinary Laguerre polynomials ($r=1$),
this result was conjectured a few years ago
by Sylvie Corteel and myself \cite{Corteel-Sokal_conj_Laguerre}
and was proven very recently by
Alex Dyachenko, Mathias P\'etr\'eolle and myself \cite{latpath_laguerre}.
Our proof is based on constructing a quadridiagonal production matrix
for the monic unsigned Laguerre polynomials $\scrlna(x)$
and then proving its total positivity;
this construction is strongly motivated by the work of
Coussement and Van Assche \cite{Coussement_03}
on the multiple orthogonal polynomials associated to weights
based on modified Bessel functions of the first kind
[cf.\ \reff{eq.laguerre.stieltjes}].
We have not yet succeeded in extending this proof to $r > 1$.

\section*{Acknowledgments}

I wish to thank the organizers of the
15th International Symposium on Orthogonal Polynomials,
Special Functions and Applications
(Hagenberg, Austria, 22--26 July 2019)
for inviting me to give a talk there;
this allowed me to meet Walter Van Assche and to discover
an unexpected connection \cite{Sokal_multiple_OP}
between branched continued fractions and multiple orthogonal polynomials,
which formed part of the motivation for this work.

I also wish to thank Kathy Driver for drawing my attention to
Hille's paper \cite{Hille_29},
and Alex Dyachenko for helpful discussions.

This research was supported in part by
the U.K.~Engineering and Physical Sciences Research Council grant EP/N025636/1.

\end{document}